\documentclass{article}
\usepackage{arxiv}

\usepackage[utf8]{inputenc}
\usepackage[T1]{fontenc}
\usepackage{microtype}

\usepackage{amsmath,amsfonts,amssymb}

\usepackage{algorithmic}
\usepackage{algorithm}

\usepackage{array}
\usepackage{booktabs}
\usepackage{makecell}
\usepackage{graphicx}
\usepackage{subfig}
\usepackage{capt-of}

\usepackage{textcomp}
\usepackage{url}
\usepackage{verbatim}
\usepackage{xcolor}

\usepackage[numbers,square,sort&compress]{natbib}
\usepackage{doi}
\usepackage{hyperref}
\hypersetup{hidelinks}

\newcommand{\cG}{\mathcal{G}}
\newcommand{\cL}{\mathcal{L}}
\newcommand{\cB}{\mathcal{B}}
\newcommand{\cN}{\mathcal{N}}
\newcommand{\cT}{\mathcal{T}}

\title{Strengthening AC Optimal Power Flow Relaxation through
Offline Bound Tightening and Linear Approximations}

\author{
    Xiangxin An \quad Constance Crozier \quad Santanu S. Dey \\
    H. Milton Stewart School of Industrial and Systems Engineering \\
    Georgia Institute of Technology \\
    Atlanta, GA 30332, USA
}

\date{}

\renewcommand{\shorttitle}{Strengthening AC Optimal Power Flow Relaxation through Offline Bound Tightening and Linear Approximations}

\hypersetup{
    pdftitle={Strengthening AC Optimal Power Flow Relaxation through Offline Bound Tightening and Linear Approximations},
    pdfauthor={Xiangxin An, Constance Crozier, Santanu S. Dey},
    pdfsubject={AC optimal power flow and convex relaxations},
    pdfkeywords={AC optimal power flow, convex relaxation, bound tightening, valid inequality}
}
\renewcommand{\undertitle}{}
\renewcommand{\headeright}{}
\begin{document}

\maketitle

\begin{abstract}
The alternating-current optimal power flow (ACOPF) problem is
fundamental to power system operation and planning, but its
nonconvexity makes globally optimal solutions of large-scale
instances computationally challenging. Convex relaxations provide
tractable lower bounds for assessing solution quality and supporting
global optimization, although their effectiveness depends strongly
on relaxation strength. This paper proposes an offline convexification
and bound tightening framework for constructing strong linear
programming (LP) relaxations for families of related ACOPF instances.
Rather than repeating the tightening procedure for each instance,
we compute variable bounds and second-order cone outer-approximation
(SOC-OA) cuts offline using a demand-slack model and reuse them for
subsequent instances with the same topology and admissible demand
variations. Computational results on PGLib instances show that the
fully offline approach achieves an average gap of $4.32\%$, compared
with $3.80\%$ for the fully online approach, while achieving an
average speedup of approximately $272\times$. The results demonstrate
that offline convexification and bound tightening can recover most
of the relaxation improvement of online tightening while significantly
reducing the computational burden for repeated ACOPF evaluations.
\end{abstract}


\section*{Nomenclature}
\subsection*{Sets and Indices}
\begin{tabular}{>{\itshape}p{0.17\columnwidth} p{0.76\columnwidth}}
$\cB$ & Set of buses $i$. \\
$\cG$ & Set of generators $k$. \\
$\cL$ & Set of transmission lines $\ell$.\\
$\cT$ & Time-period. \\
% $\cG_i$ & Set of generators located at bus $i$.\\
$\delta^+(i)$ & Set of lines whose from-side bus is $i$.\\
$\delta^-(i)$ & Set of lines whose to-side bus is $i$.\\
\end{tabular}

\subsection*{Decision Variables}
\begin{tabular}{>{\itshape}p{0.1\columnwidth} p{0.82\columnwidth}}
$p^g_{k}$, $q^g_{k}$ & Active and reactive power of generator $k$.\\
$\theta_{i}$ & Voltage phase angle at bus $i$.\\
$\theta_{ji}$ & Voltage angle difference across line $(\ell,i,j)$.\\
$V_i$ & Voltage magnitude at bus $i$.\\
$c_{ii}$ & Squared voltage magnitude at bus $i$, $c_{ii}=|V_{i}|^2$.\\
$c_{ij}$ & Voltage product, $c_{ij}=|V_{i}||V_{j}|\cos(\theta_{i}-\theta_{j})$.\\
$s_{ij}$ & Voltage product, $s_{ij}=|V_{i}||V_{j}|\sin(\theta_{i}-\theta_{j})$.\\
\end{tabular}

\begin{tabular}{>{\itshape}p{0.17\columnwidth} p{0.76\columnwidth}}
$p^{ft}_{\ell}$, $q^{ft}_{\ell}$ & Active and reactive power flow on line $\ell$ from bus $i$ to bus $j$.\\
$p^{tf}_{\ell}$, $q^{tf}_{\ell}$ & Active and reactive power flow on line $\ell$ from bus $j$ to bus $i$.\\
\end{tabular}

\subsection*{Parameters}
\begin{tabular}{>{\itshape}p{0.17\columnwidth} p{0.76\columnwidth}}
$C_k^g(\cdot)$ & Generation cost function of generator $k$.\\
$g_i, b_i$ & Shunt conductance and susceptance at bus $i$.\\
$G_\ell^{ff}, B_\ell^{ff}$ & Real and imaginary parts of the from-side self-admittance of line $\ell$.\\
$G_\ell^{ft}, B_\ell^{ft}$ & Real and imaginary parts of the from-to mutual admittance of line $\ell$.\\
$G_\ell^{tf}, B_\ell^{tf}$ & Real and imaginary parts of the to-from mutual admittance of line $\ell$.\\
$G_\ell^{tt}, B_\ell^{tt}$ & Real and imaginary parts of the to-side self-admittance of line $\ell$.\\
$p^d_{i}, q^d_{i}$ & Active and reactive power demand at bus $i$.\\
$\underline{p}^g_{k}, \overline{p}^g_{k}$ & Lower and upper active power limits of generator $k$.\\
$\underline{q}^g_{k}, \overline{q}^g_{k}$ & Lower and upper reactive power limits of generator $k$.\\
$\underline{RD}_k, \overline{RU}_k$ & Ramp-down and ramp-up limits of generator $k$.\\
$\underline{V}_i, \overline{V}_i$ & Lower and upper voltage magnitude limits at bus $i$.\\
$\overline{S}_{\ell}$ & Thermal limit of line $\ell$.
\end{tabular}

\section{Introduction}

The alternating-current optimal power flow (ACOPF) problem is one of the fundamental optimization problems in power system operation and planning. It seeks to minimize generation costs while satisfying network physics, operational limits, and real and reactive power balance constraints. Due to the nonlinearity induced by Ohm's and Kirchhoff's laws, ACOPF is a nonlinear and nonconvex optimization problem. Consequently, solving ACOPF to global optimality is computationally challenging, especially for large-scale networks. To address this challenge, convex relaxations have received significant attention as a means of producing tractable lower bounds and certifying solution quality, with widely studied approaches including linear programming (LP) relaxations~\cite{coffrin2014linear}, second-order cone programming (SOCP) relaxations~\cite{1664986,7056568}, quadratic convex (QC) relaxations~\cite{7271127}, and semidefinite programming (SDP) relaxations~\cite{bai2008semidefinite}. A comprehensive survey of relaxations and approximations of the power flow equations can be found in~\cite{molzahn2019survey}.

For a minimization problem such as ACOPF, a convex relaxation provides a lower bound on the globally optimal objective value, while any AC-feasible solution provides an upper bound. The difference between these bounds, commonly referred to as the relaxation gap, indicates the quality of the relaxation: a smaller gap implies that the relaxation more closely approximates the original nonconvex problem. However, the strength of convex relaxations can vary significantly across instances. To improve the strength of these relaxations, several valid inequalities, such as arctangent envelopes~\cite{kocuk2016strong} and lifted nonlinear cuts~\cite{10707070}, have been proposed to further strengthen the relaxation by adding more constraints to the model. At the same time, loose variable bounds can lead to large relaxation gaps. Therefore, bound tightening has become a powerful strengthening tool, since tighter variable bounds can directly improve the quality of convex envelopes and reduce the feasible region of the relaxation.

Bound tightening methods for ACOPF can be broadly divided into two categories: analytical bound tightening and optimization-based bound tightening. Analytical bound-tightening methods exploit the specific structure of AC power-flow equations to derive explicit or nearly closed-form tightening rules. This type of method is usually applied on voltage magnitudes, phase-angle differences, nodal injections, and line flows in ACOPF problems on different relaxations such as LP, QC, and SDP~\cite{7328765,coffrin2015strengthening,7763860}. 

Optimization-based bound tightening has also been widely applied and is one of the most effective procedures for reducing variable domains \cite{6814337,8500764,bynum2018strengthened,10384116}. The standard computational bound tightening procedure solves auxiliary optimization problems to obtain the tightest bounds that are valid for all solutions in the relaxation, by minimizing and maximizing each selected variable over the relaxed feasible region. However, these bounds are generally computed with respect to the feasible region of a particular operating instance and therefore must typically be recomputed when the operating conditions change. This procedure can be further strengthened by adding an objective-cutoff constraint to cut-off all solutions whose objective values are worse than the best-known feasible upper bound. This variant is commonly referred to as objective-based bound tightening. To avoid confusion with the broader term optimization-based bound tightening, we hereby avoid using the abbreviation OBBT, and refer to the standard version as feasibility-based bound tightening (FBBT), since it only uses the relaxed feasible region and does not exploit information from the objective function. 

Previous studies have shown that FBBT can substantially tighten voltage- and phase-related variable domains and improve the strength of LP, SOCP, and QC relaxations for AC power-flow and ACOPF problems~\cite{6814337,8500764,bynum2018strengthened,10384116}. Tighter bounds can also improve the performance of mixed-integer AC power-sysetem optimization problems, such as AC unit commitment (ACUC) and AC optimal transmission switching (ACOTS), by strengthening relaxations and facilitating branch-and-bound pruning, as demonstrated in applications such as optimal transmission switching and generator switching~\cite{guo2025tightening,aigner2023solving}. 

Despite its effectiveness, FBBT remains computationally expensive as it requires solving two auxiliary optimization problems for each selected variable, and these problems may need to be solved repeatedly over several iterations. Consequently, the total computational cost can become prohibitive for large-scale networks. Existing studies have tried to improve the scalability through network decomposition~\cite{bynum2021decomposition}, selective variable tightening~\cite{shao2025efficient, CENGIL2022108275}, and parallel computation~\cite{8667386}. These approaches substantially reduce the cost of performing bound tightening on an individual instance, but they do not eliminate the need to repeat the tightening procedure across multiple related load instances. 

This limitation becomes particularly important in applications such as day-ahead operation, contingency analysis, scenario-based planning, and multi-period optimization, where many similar ACOPF instances must be solved repeatedly within a limited time. Although these instances are often closely related, conventional bound-tightening approaches generally treat each one independently, leading to repeated solving similar optimization problems and substantial redundant computation. 

Related work has explored offline preprocessing over a range of operating conditions. In particular, Aquino et al.~\cite{9654442} introduced a demand-flexibility formulation to identify redundant ACOPF constraints across varying load conditions, allowing the resulting reduced formulation to be reused for subsequent instances. While their approach focuses on constraint screening and model reduction, our work extends the use of demand flexibility to offline relaxation strengthening through reusable variable bounds and SOC outer-approximation cuts.

Motivated by the need for efficient relaxation strengthening across
related instances, we propose an offline convexification and bound-tightening framework that exploits similarities across related ACOPF instances. The main contributions of this paper are as follows:
\begin{enumerate}
    \item We develop an unified demand-slack-based offline bound tightening approach for a family of instances sharing the same network topology and varying demands. Computational results show that the resulting bounds retain relaxation strength comparable to those generated from the nominal model.
    \item We propose an integrated convexification and bound tightening framework to construct a stronger LP relaxation of ACOPF. FBBT is used to reduce the domains of the lifted voltage-related variables, and the arctangent envelopes and SOC-OA cuts, building on \cite{kocuk2016strong, bienstock2026accurate}, further closes the gap between the LP and SOCP relaxations without sacrificing the linear structure.
    \item Instead of performing expensive strengthening procedures independently for each new demand realization, we apply the proposed framework fully offline on the demand slack model to get reusable variable bounds and cuts over a prescribed range of demand variations. 
\end{enumerate}

The main advantage of the proposed framework is that it separates the expensive strengthening phase from the online solution phase. Since the generated bounds and cuts are constructed to remain globally valid for the targeted family of instances, they can be reused across multiple operating scenarios without sacrificing relaxation validity. Furthermore, because the framework relies on network-level structure rather than a specific single-period formulation, it can be naturally extended to mixed-integer variants of ACOPF, including AC unit commitment and optimal transmission switching.

\section{Problem formulation}
Consider a transmission-level power network $\cN = (\cB,\cL)$, where $\mathcal{B}$ denotes the set of buses and $\mathcal{L}$ denotes the set of transmission lines. The active and reactive power flows on each transmission line are nonlinear functions of bus voltage magnitudes $V_i$, $V_j$ and voltage angle differences $\theta_{ij}$. To obtain a compact formulation of the flow power, we introduce the following lifted variables:
\begin{subequations}\label{eq:cs_expression}
\begin{align}
c_{ii}&=|V_{i}|^2, & \forall i \in \cB\label{eq:1A}\\
c_{ij} &= |V_i||V_j| \cos (\theta_{i} - \theta_{j}), & \forall (\ell, i, j) \in \cL \label{eq:1B}\\
s_{ij} &= -|V_i||V_j| \sin (\theta_{i} - \theta_{j}), & \forall (\ell, i, j) \in \cL \label{eq:1C}
\end{align}
\end{subequations}

\subsection{Exact model}
Using the lifted variables in \eqref{eq:cs_expression}, the single-period ACOPF problem can be formulated in \eqref{model:ACOPF_single_exact}. The objective function \eqref{eq:model_obj_opf} minimizes the total generation cost, where $C_k^g(p_{k}^g)$ represents the cost function of generator $k$, typically modeled as a quadratic function of active power generation $p_{k}^g$. Constraints \eqref{eq:jabr_p_fr}-\eqref{eq:jabr_q_to} represent the active and reactive power flow representations in both directions of each transmission line in the branch-flow model. \eqref{eq:jabr_p_balance}-\eqref{eq:jabr_q_balance} enforce active and reactive power balance at each bus. For each generator, the active and reactive power outputs are restricted by their operating physical limits in \eqref{eq:gen_power_limit}, while the bus voltage magnitude is also bounded in \eqref{eq:jabr_vm_limit}. The thermal limits \eqref{eq:thermal_limit} ensure that apparent power flows remain within limits to prevent overheating and instability of the line. Finally, by introducing the $c,s$ variables, \eqref{eq:cs_rotated_ring} and \eqref{eq:jabr_tan_rel} enforce the exact nonlinear relationships among the lifted variables and voltage angle differences.

\begin{subequations}\label{model:ACOPF_single_exact}
    \begin{align}  
    &\min\;\sum_{k\in\mathcal{G}}C_k^g\left(p^{g}_{k}\right) \label{eq:model_obj_opf}\\[2pt]
      &\forall \text{ transmission line } (\ell,i,j)\in\mathcal{L}: \notag\\
      &p^{ft}_{\ell} = G^{ff}_{\ell}\, c_{ii} + G^{ft}_{\ell}c_{ij} - B^{ft}_{\ell}\,s_{ij} \label{eq:jabr_p_fr}\\
      &q^{ft}_{\ell} = -B^{ff}_{\ell}\,c_{ii} - B^{ft}_{\ell}c_{ij} - G^{ft}_{\ell}s_{ij} \label{eq:jabr_q_fr}\\
      &p^{tf}_{\ell} = G^{tt}_{\ell}\,c_{jj} + G^{tf}_{\ell}c_{ij} + B^{tf}_{\ell}s_{ij} \label{eq:jabr_p_to}\\
      &q^{tf}_{\ell} = -\,B^{tt}_{\ell}\,c_{jj} - B^{tf}_{\ell}c_{ij} + G^{tf}_{\ell}\,s_{ij} \label{eq:jabr_q_to}\\
      &\forall \text{ bus } i \in \mathcal{B}: \notag \\
      &\sum_{k \in \mathcal{G}_i} p^{g}_{k}-p^{d}_{i} = g_{i}\,c_{ii}
      +\sum_{\ell\in\delta^+(i)}p^{ft}_{\ell}
      +\sum_{\ell\in \delta^-(i)}p^{tf}_{\ell} \label{eq:jabr_p_balance}\\
      &\sum_{k \in \mathcal{G}_i} q^{g}_{k}-q^{d}_{i} = -\,b_{i}\,c_{ii}
      +\sum_{\ell\in\delta^+(i)}q^{ft}_{\ell}
      +\sum_{\ell\in\delta^-(i)}q^{tf}_{\ell} \label{eq:jabr_q_balance} \\
      & \forall \text{ generator } k \in \mathcal{G}:  \underline{p}^{g}_{k} \le p^{g}_{k}\le\overline{p}^{g}_{k}, \quad \underline{q}^{g}_{k} \le q^{g}_{k}\le\overline{q}^{g}_{k} \label{eq:gen_power_limit}\\
      & \forall \text{ bus } i \in \mathcal{B}: \underline{V}_{i}^2 \le c_{ii} \le\overline{V}_{i}^2 \label{eq:jabr_vm_limit}\\
      &\forall \text{ transmission line } (\ell,i,j)\in\mathcal{L}: \notag\\
      & (p^{ft}_{\ell})^{2}+(q^{ft}_{\ell})^{2} \le\overline{S}_{\ell}^{2}, \quad (p^{tf}_{\ell})^{2}+(q^{tf}_{\ell})^{2} \le\overline{S}_{\ell}^{2} \label{eq:thermal_limit}\\
      & c_{ij}^{2}+s_{ij}^{2} = c_{ii}\,c_{jj} \label{eq:cs_rotated_ring}\\
      & \theta_{ji} = \theta_{j} - \theta_{i} = \arctan \frac{s_{ij}}{c_{ij}} \label{eq:jabr_tan_rel}
    \end{align}
\end{subequations}

\subsection{SOCP relaxation}

Following the classical SOCP relaxation of Jabr \cite{1664986}, we relax the nonlinear equality \eqref{eq:cs_rotated_ring} into a rotated second-order cone constraint and omit the voltage angle difference constraint \eqref{eq:jabr_tan_rel}. The resulting SOCP relaxation is given by:
\begin{subequations}\label{model:SOCP}
\begin{align}
    & \eqref{eq:model_obj_opf}-\eqref{eq:thermal_limit} \notag\\
    & c_{ij}^{2}+s_{ij}^{2} \le c_{ii}\,c_{jj} \label{eq:cs_rotated_cone}
\end{align}    
\end{subequations}
This relaxation preserves the linear power flow equations, nodal balance constraints, generator limits, voltage magnitude bounds, and thermal constraints, while replacing the nonconvex equality with a convex conic inequality.

\subsection{LP relaxation}

An LP relaxation can be obtained by replacing quadratic and bilinear terms in the exact ACOPF formulation \eqref{model:ACOPF_single_exact} with their McCormick envelope approximations. Specifically, for a bilinear term $w=xy$, where $\underline{x}\le x\le \overline{x}$ and $\underline{y}\le y\le \overline{y}$, the McCormick relaxation $M(w=xy)$is defined by the following inequalities:
\begin{subequations}
    \begin{align}
    &w \ge x\underline{y} + \underline{x}y - \underline{x}\underline{y}, &w \ge \overline{x}y + x\overline{y} - \overline{x}\overline{y} \label{eq:McCormick1}\\
    & w \le x\overline{y} + \underline{x}y - \underline{x}\overline{y}, & w \le  x\underline{y} + \overline{x}y - \overline{x}\underline{y} \label{eq:McCormick4}
    \end{align}
\end{subequations}
Using these envelopes, the quadratic terms in the thermal limits and the bilinear term in \eqref{eq:cs_rotated_ring} are replaced by auxiliary variables. The resulting LP relaxation can be written as:
\begin{subequations}\label{model:McCormick}
    \begin{align}
        & \eqref{eq:model_obj_opf}-\eqref{eq:jabr_vm_limit} \notag \\
        &\forall (\ell, i, j) \in \cL: \notag \\
        & P_{\ell}^{ft} + Q_{\ell}^{ft} \le \overline{S}_{\ell}^2, P_{\ell}^{tf} + Q_{\ell}^{tf} \le \overline{S}_\ell^2,  \label{eq:LP_thermal_limit}\\
        &C_{ij} +S_{ij} = D_{ij}, \label{eq:CS_relationship} \\
        & M(C_{ij} = c_{ij}^2), M(S_{ij} = s_{ij}^2), M(D_{ij} = c_{ii}c_{jj}) \label{eq:mccormick_cs}\\
        & M(P_{\ell}^{ft}=(p_{\ell}^{ft})^2), M(Q_{\ell}^{ft}=(q_{\ell}^{ft})^2) \label{eq:mccormick_fr} \\
        & M(P_{\ell}^{tf}=(p_{\ell}^{tf})^2), M(Q_{\ell}^{tf}=(q_{\ell}^{tf})^2) \label{eq:mccormick_to}
    \end{align}
\end{subequations}
Here, $M(\cdot)$ denotes the McCormick envelope defined by \eqref{eq:McCormick1}-\eqref{eq:McCormick4}. This LP relaxation provides a computationally efficient approximation of the ACOPF feasible region and is particularly useful when combined with bound-tightening procedures, since tighter variable bounds directly strengthen the McCormick envelopes.

\section{Offline convexification and bound tightening Framework}

The proposed offline convexification and bound tightening framework exploits the observation that many ACOPF instances share the same network topology and differ mainly in their demand profiles and generator cost functions. Instead of repeatedly performing bound tightening and cut generation for each new instance, these computationally expensive steps are carried out once offline over a prescribed demand range. The resulting variable bounds and valid cuts are then stored and reused to strengthen future instances with the same topology.

\begin{figure}[!h]
    \centering
    \includegraphics[width=0.6\linewidth]{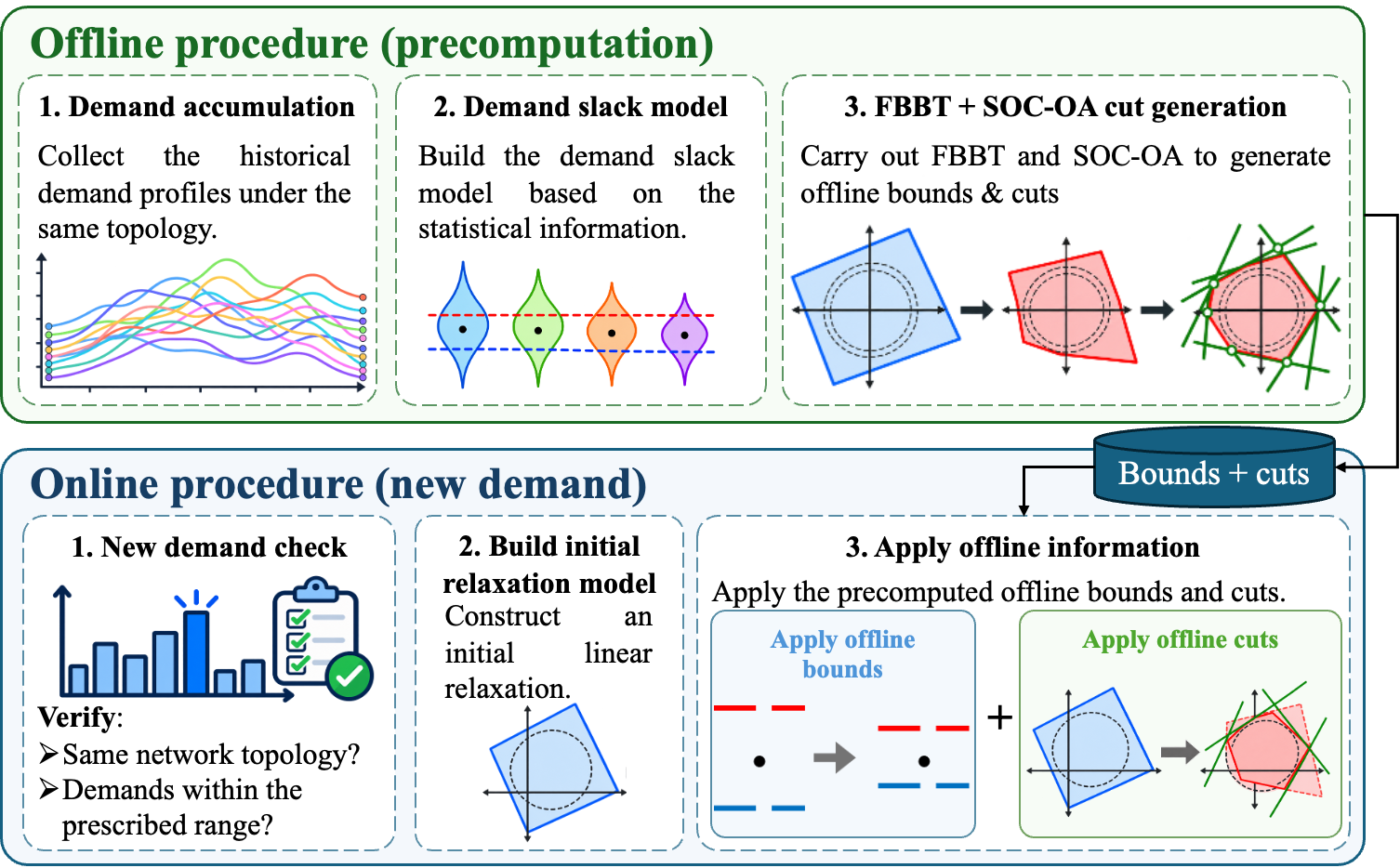}
    \caption{The offline convexification and bound tightening framework.}
    \label{fig:framework}
\end{figure}

Figure~\ref{fig:framework} illustrates the overall workflow of the proposed framework. The upper panel shows the offline precomputation procedure, where historical demand profiles are aggregated to construct a demand slack model. Bound tightening and SOC outer approximation (SOC-OA) are then performed on this model to generate reusable offline bounds and cuts. Specifically, FBBT iteratively computes and tightens the variable bounds, while SOC-OA generates linear inequalities tangent to the SOC constraints. The lower panel shows the online procedure. Given a new demand instance with the same topology and demands within the prescribed offline range, an initial linear relaxation is constructed and strengthened directly using the precomputed bounds and cuts, completely bypassing the expensive online bound-tightening and cut-generation procedures.

\subsection{Linear approximation enhancement tools}
We first introduce two standard classes of linear approximation tools used in the proposed framework. These tools play different roles in different phases of the algorithm. Arctangent envelopes are used during the offline FBBT iterations to strengthen the relaxation and improve the quality of the computed bounds. They can also be used in any other relaxations. SOC outer-approximation cuts, on the other hand, are generated after the FBBT phase. Since they are derived from conic constraint geometry, these cuts can be computed offline and reused across future instances.

\subsubsection{Arctangent envelopes (AT)}
The arctangent envelopes are built to outer approximate the constraints \eqref{eq:jabr_tan_rel} using methods from \cite{kocuk2016strong}. In our previous paper, this type of inequality provides improvements on the LP relaxation \cite{11272395}, so we continue using it to enhance the LP and SOCP relaxation. Since $c_{ij}$ and $s_{ij}$ are bounded, the region in 3D that we need to describe is bounded by a 2D box in the $c,s$ plane.
\begin{align}
   \mathcal{A}_{ij} := \{ & (c_{ij},s_{ij},\theta_{ji}) \in \mathbf{R}^3: \theta_{ji} = \arctan(s_{ij}/c_{ij}), \notag\\
    &(c_{ij}, s_{ij}) \in [\underline{c}_{ij},\overline{c}_{ij}] \times [\underline{s}_{ij}, \overline{s}_{ij}]\}, \notag 
\end{align}
To obtain the convex hull of $\mathcal{A}_{ij}$, we project the four corner points of $\mathcal{B}_{ij}$ in 2D onto the arctangent surface $\theta_{ji} = \arctan(s_{ij}/c_{ij})$, then we obtain 4 corner points of $\mathcal{A}_{ij}$: $z_{ij}^1 = (\underline{c}_{ij}, \overline{s}_{ij}, \arctan(\overline{s}_{ij}/\underline{c}_{ij}))$, $z_{ij}^2 = (\overline{c}_{ij},\overline{s}_{ij}, \arctan(\overline{s}_{ij}/\overline{c}_{ij}))$, $z_{ij}^3 = (\overline{c}_{ij}, \underline{s}_{ij}, \arctan(\underline{s}_{ij}/\overline{c}_{ij}))$, and $z_{ij}^4 = (\underline{c}_{ij}, \underline{s}_{ij}, \arctan(\underline{s}_{ij}/\underline{c}_{ij}))$.
The upper envelopes of $\mathcal{A_{ij}}$ can be defined by two supporting planes that pass through the point triples $\{z_{ij}^1, z_{ij}^2, z_{ij}^3\}$ and $\{z_{ij}^1, z_{ij}^3, z_{ij}^4\}$, denoted as $\theta_{ji}^1 = \gamma_1 + \alpha_1 c_{ij}+ \beta_1 s_{ij}$ and  $\theta_{ji}^2 = \gamma_2 + \alpha_2 c_{ij} + \beta_2 s_{ij}$: 
\vspace{0 pt}
\begin{equation}
\begin{split}
    & \gamma_k' + \alpha_k c_{ij} + \beta_k s_{ij} \ge \theta_{ji},\\
    & \text{where } \gamma_k' = \gamma_k + \Delta \gamma_k, \\
    & \Delta \gamma_k =\max_{(c_{ij},s_{ij})\in \mathcal{B}_{ij}}\{\arctan(s_{ij}/c_{ij})-\theta_{ji}^k\},\\
    & \theta_{ji}^k = \gamma_k + \alpha_k c_{ij} + \beta_k s_{ij}, \text{and } k = 1,2.
\end{split} \label{eq:arctan_upper}
\end{equation}

Similarly, the lower envelope is constructed by fitting the planes $\theta_{ji}^3 = \gamma_3 + \alpha_3 c_{ij} + \beta_3 s_{ij}$ and  $\theta_{ji}^4 = \gamma_4 + \alpha_4 c_{ij} + \beta_4 s_{ij}$ through the triples $\{z_{ij}^1, z_{ij}^2, z_{ij}^4\}$ and $\{z_{ij}^2, z_{ij}^3, z_{ij}^4\}$.
\vspace{0 pt}
\begin{equation}
\begin{split}
    & \gamma_k' + \alpha_k c_{ij} + \beta_k s_{ij} \le \theta_{ji},\\
    & \text{where } \gamma_k' = \gamma_k - \Delta \gamma_k,\\
    &\Delta \gamma_k= \max_{(c_{ij},s_{ij})\in \mathcal{B}_{ij}}\{\theta_{ji}^k-\arctan(s_{ij}/c_{ij})\},\\
    & \theta_{ji}^k = \gamma_k + \alpha_k c_{ij} + \beta_k s_{ij}, \text{and } k = 3,4.
\end{split} \label{eq:arctan_lower}
\end{equation}

Noticeable, the arctangent envelopes are independent of the load demands, and only depend on the bounds of $c,s$ variables.

\subsubsection{SOC outer-approximation (SOC-OA)}
The second class of linear approximations is based on outer approximations of SOC constraints. After offline bound tightening, the resulting model is still an LP relaxation if the conic constraints are not imposed explicitly. To further close the gap between the LP relaxation and the SOCP relaxation while preserving a linear formulation, we extend the cutting-plane approach in \cite{bienstock2026accurate} to approximate SOC constraints with linear inequalities. 

For the standard second-order cone $\mathcal{S}_r:=\{ (x, y) \in \mathbb{R}^2: x^2+y^2 \le r^2\}$, given a current solution $(x^*, y^*) \notin \mathcal{S}_r$, the corresponding separating hyperplane is defined as: 
\begin{equation}
    x^* x + y^*y \le r\sqrt{(x^*)^2+(y^*)^2} \label{eq:soc_cutting_plane}
\end{equation}

For the rotated second-order cone $\mathcal{C}:=\{ (x, y, w, z) \in \mathbb{R}^2 \times\mathbb{R}^2_+: x^2+y^2 \le wz\}$ with $(x^*, y^*, w^*, z^*) \notin \mathcal{C}$, $w\ge 0, z \ge 0 $. The hyperplane is defined as: 
\begin{equation}
    4x^* x + 4y^*y + ((w^* - z^*) - n_0)w +(-(w^* - z^*) - n_0)z \le 0 \label{eq:rotated_soc_cutting_plane}
\end{equation}
where $n_0 := ||(2x^*, 2y^*, w^* - z^*) ||_2$.
Starting from a simple LP relaxation, the cutting-plane algorithm iteratively solves the current relaxation, identifies violated conic constraints, and adds the corresponding linear cuts. With appropriate cut-management strategies, this procedure yields a progressively tighter LP approximation of the SOCP relaxation. Moreover, since these cuts are derived from the conic constraint geometry rather than from the specific load realization, they remain valid across instances as long as the associated transmission line remains in service.

\subsection{Offline bound tightening}
\subsubsection{Base algorithm}

Bound tightening is an effective way to strengthen convex relaxations of ACOPF, especially when McCormick envelopes or other variable-bound-dependent approximations are used. Classical FBBT is typically performed on voltage magnitude and voltage angle variables. In this paper, we focus on tightening the lifted variables $c_{ij}$ and $s_{ij}$, since these variables directly determine the strength of the McCormick envelopes and arctangent envelopes.

Given a relaxation $\mathcal{S}$ of the ACOPF feasible region, the standard FBBT tightens the bounds of each selected variable by solving two auxiliary optimization problems over the current relaxation. For each selected variable $x\in\{c_{ij},s_{ij}: \forall (\ell, i,j) \in \mathcal{L}\}$, its candidate lower and upper bounds are obtained as
\[
\underline{x}'=\min_{\mathcal{S}} x, \qquad \overline{x}'=\max_{\mathcal{S}} x.
\] 
The bounds are then updated by intersecting these values with the current variable domain,
\[
\underline{x}^{k} = \max\{\underline{x}',\underline{x}^{k-1}\},
\qquad
\overline{x}^{k} = \min\{\overline{x}',\overline{x}^{k-1}\},
\]
where $k$ denotes the bound-tightening iteration. The relaxation $\mathcal{S}$ is then updated with the tightened bounds, and the process is repeated until the bound changes fall below a prescribed tolerance or a maximum number of iterations is reached. 

The main drawback of FBBT is its computational cost. Since two auxiliary optimization problems must be solved for each selected variable, the number of optimization problems grows with both the network size and the number of tightening iterations. Moreover, when operating conditions, such as demand profile, change, the procedure must be performed again from scratch for the new instance. This limits its practical applicability in real-time or repeated operational settings. Therefore, repeating these largely similar optimization problems across a family of related instances can therefore result in substantial redundant computation, motivating the demand slack model and the offline bound-tightening framework developed in the following sections.

\subsubsection{Demand slack model}
In many power system applications, the network topology, voltage limits, generator capacities, and line parameters remain unchanged over a large number of operating instances. The main variations across instances typically come from load demand profiles and generation cost functions. This observation motivates the construction of an offline demand slack model, which captures a representative range of operating conditions and enables reusable bound-tightening computations.

Consider a set of historical ACOPF instances with identical topology. For each bus $i$, we construct a nominal demand value and a symmetric admissible demand range from the historical observations. Let:
\begin{equation}\label{eq:demand range}
\begin{aligned}
p_{i,\min}^d
    := \min_{r\in\mathcal I} p_i^{d,(r)},\quad
&
p_{i,\max}^d
    := \max_{r\in\mathcal I} p_i^{d,(r)},\\
q_{i,\min}^d
    := \min_{r\in\mathcal I} q_i^{d,(r)},\quad
&
q_{i,\max}^d
    := \max_{r\in\mathcal I} q_i^{d,(r)},\\
\tilde p_i^d
    := \frac{p_{i,\min}^d+p_{i,\max}^d}{2},\quad
&
\tilde q_i^d
    := \frac{q_{i,\min}^d+q_{i,\max}^d}{2},\\
\alpha_i^p
    := \frac{p_{i,\max}^d-p_{i,\min}^d}
              {p_{i,\max}^d+p_{i,\min}^d},\quad
&
\alpha_i^q
    := \frac{q_{i,\max}^d-q_{i,\min}^d}
              {q_{i,\max}^d+q_{i,\min}^d}.
\end{aligned}
\end{equation}
where $\alpha_i^p$ and $\alpha_i^q$ measure the relative half-widths of
the historical active and reactive demand ranges, respectively. The demand slack model is built from the same single-period ACOPF relaxation by changing the node balance from \eqref{eq:jabr_p_balance}\eqref{eq:jabr_q_balance} to \eqref{eq:jabr_p_slack}\eqref{eq:jabr_q_slack}: 
\begin{subequations}\label{model:agg}
    \begin{align}
    & \text{Any relaxations of ACOPF, i.e. \eqref{model:SOCP} or \eqref{model:McCormick} without \eqref{eq:jabr_p_balance}\eqref{eq:jabr_q_balance}} \notag \\
    &\sum_{k \in \mathcal{G}_i} p^{g}_{k}-\tilde{p}^{d}_{i} + \gamma_i^{p} = g_{i}\,c_{ii}
      +\sum_{\ell\in\delta^+(i)}p^{ft}_{\ell}
      +\sum_{\ell\in \delta^-(i)}p^{tf}_{\ell} \label{eq:jabr_p_slack}\\
      &\sum_{k \in \mathcal{G}_i} q^{g}_{k}-\tilde{q}^{d}_{i} + \gamma_i^{q} = -\,b_{i}\,c_{ii}
      +\sum_{\ell\in\delta^+(i)}q^{ft}_{\ell}
      +\sum_{\ell\in\delta^-(i)}q^{tf}_{\ell} \label{eq:jabr_q_slack}\\
      & -\alpha_i^p\tilde{p}_i^d \le \gamma_i^{p} \le \alpha_i^p \tilde{p}_i^d, \quad -\alpha_i^q \tilde{q}_i^d \le \gamma_i^{q}\le \alpha_i^q \tilde{q}_i^d.\label{eq:pq_slack_bounds}
\end{align}
\end{subequations}
Here, $\gamma_i^{p}$ and $\gamma_i^{q}$ represent the deviations from the nominal active and reactive demand, respectively. $\alpha_i^p$ and $\alpha_i^q$ are estimated from historical demand data. With these slack variables, the demand slack model contains all demand realizations whose loads lie within the specified ranges.

\subsubsection{Offline bound tightening algorithm}
After constructing the offline demand slack model, we perform FBBT on that model to compute tightened bounds for the selected $c_{ij}$ and $s_{ij}$ variables. Since the slack model explicitly incorporates admissible load perturbations, the resulting bounds are valid not only for the median-demand instance but also for any future instance whose demand profile lies within the prescribed slack ranges. Therefore, in the online procedure, these offline-computed bounds can be directly reused without rerunning the full FBBT procedure.

It is also important to note that generation cost functions may vary across operating instances. To ensure that the offline-computed bounds remain representative across different cost realizations, we do not impose objective cutoff constraints in the offline FBBT procedure. The resulting offline bound tightening framework is summarized in Algorithm~\ref{alg:offline}.

\begin{algorithm}[H]
\caption{Offline Bound Tightening Framework}
\label{alg:offline}
\begin{algorithmic}[1]

\REQUIRE
Historical instances $\mathcal{I}$ with the same topology; maximum iterations $K$; tolerance \texttt{bound\_tol}, \texttt{improve\_tol}; arctangent-update set $\mathcal{K}_{\rm atan}$.

\ENSURE
Offline tightened bounds for selected $(c,s)$ variables.

\STATE \textbf{Base model construction:}
Build base ACOPF relaxation model \eqref{model:McCormick} under the given network topology.

\STATE \textbf{Demand statistics:}
For each bus $i \in \mathcal{B}$, compute the demand range parameters in \eqref{eq:demand range}: $\tilde{p}_i^d, \alpha_i^p; \tilde{q}_i^d, \alpha_i^q$.

\STATE \textbf{Slack model:}
Construct model $\mathcal{S}^0$ using \eqref{model:agg}.

\STATE Define the set of variables to be tightened:
\[
\mathcal{X}:=
\{c_{ij}, s_{ij}: (\ell,i,j)\in\mathcal{L}\}.
\]
\STATE Initialize bounds:
$\underline{x}^{0}:=\underline{x},
\quad
\overline{x}^{0}:=\overline{x},
\quad \forall x\in\mathcal{X}$.

\FOR{$k=1,\ldots,K$}
    \IF{$k\in\mathcal{K}_{\mathrm{atan}}$}
        \STATE Add/update arctangent envelopes \eqref{eq:arctan_upper}--\eqref{eq:arctan_lower}.
    \ENDIF

    \STATE \textbf{Parallel bound computing:}
    For each $x\in\mathcal{X}$, solve:
    \[
    \underline{x}' := \min\{x: x\in\mathcal{S}^{k-1}\},
    \quad
    \overline{x}' := \max\{x: x\in\mathcal{S}^{k-1}\}.
    \]

    \STATE \textbf{Bound update:} For each $x\in\mathcal{X}$, compute the candidate bounds:
    \[
    \underline{x}^{\mathrm{cand}}:= \max\{\underline{x}^{k-1},\underline{x}'\}, \quad
    \overline{x}^{\mathrm{cand}} := \min\{\overline{x}^{k-1},\overline{x}'\}.
    \]
        \IF{$\overline{x}^{\mathrm{cand}}
              -\underline{x}^{\mathrm{cand}}
              \geq \texttt{bound\_tol}$}
            \STATE Set
            $\underline{x}^{k}:=\underline{x}^{\mathrm{cand}}$,
            $\overline{x}^{k}:=\overline{x}^{\mathrm{cand}}$.
        \ELSE
            \STATE Retain the previous bounds:
            $\underline{x}^{k-1}$, and $\overline{x}^{k-1}$.
        \ENDIF

    \STATE \textbf{Inequality update:}
    Update McCormick envelopes using the new bounds
    $\{\underline{x}^{k},\overline{x}^{k}\}_{x\in\mathcal{X}}$.
    \STATE Compute 
    \begin{equation}
    \Delta_k:= \max_{x\in\mathcal{X}}
    \left\{
    |\underline{x}^{k}-\underline{x}^{k-1}|,
    |\overline{x}^{k}-\overline{x}^{k-1}|
    \right\}. \label{eq:delta_k}
    \end{equation}
    \IF{
    $
    \Delta_k \le \texttt{improve\_tol}$
    }
        \STATE \textbf{break}
    \ENDIF

\ENDFOR

\RETURN Tightened bounds 
$\{\underline{x}^{k},\overline{x}^{k}\}_{x\in\mathcal{X}}$.

\end{algorithmic}
\end{algorithm}

The proposed offline bound tightening framework also incorporates three computational enhancements. First, arctangent envelope cuts can be added at selected iterations to strengthen the relaxation without imposing them at every iteration. Second, instead of tightening variables sequentially, each FBBT iteration fixes the current relaxation and solves the bound-tightening subproblems for different variables in parallel. The variable bounds are updated only after all parallel subproblems in the current iteration have been completed. This avoids order-dependent tightening within one iteration and improves computational efficiency. Third, after each bound update, the McCormick envelopes are updated using the newly tightened bounds. The arctangent envelopes also depend on the bounds of the $c_{ij}$ and $s_{ij}$ variables. However, due to the additional computational cost of constructing arctangent envelopes and potential numerical stability concerns, these envelopes are not updated at every iteration. Instead, they are updated only at selected iterations specified by the set $\mathcal{K}_{\mathrm{atan}}$. 

\subsection{Offline reusable SOC cut generation}

After offline bound tightening terminates, the tightened bounds are used to construct a stronger LP relaxation. Building on the approach in~\cite{bienstock2026accurate}, we generate SOC outer-approximation (SOC-OA) cuts to further close the gap between the LP and SOCP relaxations while preserving a linear formulation.

Although SOC-OA cuts generated under a fixed nominal demand remain globally valid~\cite{bienstock2026accurate}, we decide to perform the cut-generation procedure on the demand-slack model \eqref{model:agg}. Starting from the McCormick relaxation \eqref{model:McCormick} constructed using the offline tightened bounds with the demand slack \eqref{model:agg}, we first add arctangent envelopes \eqref{eq:arctan_upper}--\eqref{eq:arctan_lower} to strengthen the LP relaxation used for cut generation. We then iteratively solve the resulting LP, identify violated SOC and rotated-SOC constraints, and add the corresponding separating hyperplanes
\eqref{eq:soc_cutting_plane}--\eqref{eq:rotated_soc_cutting_plane}.
The procedure terminates when no SOC violation exceeds the prescribed tolerance, as summarized in Algorithm~\ref{alg:offline_soc_cuts}.

\begin{algorithm}[H] 
\caption{Offline Reusable SOC Cut Generation} \label{alg:offline_soc_cuts} 
\begin{algorithmic}[1] 
\REQUIRE offline tightened bounds $\{\underline{x}^{K},\overline{x}^{K}\}_{x\in\mathcal{X}}$; SOC-OA parameters in Table~\ref{tab:offline_parameters}.
\ENSURE Reusable SOC outer-approximation cuts. 
\STATE Build demand slack model using offline tightened bounds. 
\STATE Add arctangent envelopes \eqref{eq:arctan_upper}--\eqref{eq:arctan_lower} for each selected $(c_{ij},s_{ij})$ pair using the offline tightened bounds. 
\STATE Initialize iteration counter $k:=0$.
\REPEAT 
\STATE Solve the current LP relaxation. 
\STATE Identify violated constraints of \eqref{eq:thermal_limit}\eqref{eq:cs_rotated_cone} with violation exceeding \texttt{cut\_vio} and select the most violated fractions $p_{\rm thm}$ and $p_{\rm jabr}$. 
\STATE Generate the corresponding cuts use \eqref{eq:soc_cutting_plane}\eqref{eq:rotated_soc_cutting_plane}. 
\STATE Filter nearly parallel cuts using $\varepsilon_{\rm par}$.
\STATE Add these generated hyperplanes to the LP relaxation and remove inactive cuts according to \texttt{cut\_age} and \texttt{cut\_value}.
\STATE $k\leftarrow k+1$.
\UNTIL{no SOC violation exceeds \texttt{cut\_vio}} 
\STATE Store the generated SOC outer-approximation cuts. 
\RETURN Reusable SOC cuts for future instances with the same topology. 
\end{algorithmic} 
\end{algorithm}

The generated SOC-OA cuts are stored together with the offline tightened bounds and directly applied to subsequent instances. The arctangent envelopes, however, are not stored as reusable cuts; they are regenerated online using the bounds available for each instance. This step is computationally inexpensive, requiring only the construction of linear planes from the current $c$-$s$ bounds. Consequently, the online stage avoids repeating both FBBT and the SOC cut-generation procedure.

\section{Computational Results}
\subsection{Data and parameters settings}

The test instances are selected from the PGLib library, and we focus primarily on instances with fewer than 1000 buses. For each original instance, we generate 100 perturbed load scenarios by independently sampling each active demand within $\pm 10\%$ of its nominal value, following~\cite{xavier2021learning}. Specifically,
\[
p_i^d \sim \mathcal{U}\bigl((1-10\%)\tilde p_i^d,\,(1+10\%)\tilde p_i^d\bigr),
\]
where $\tilde p_i^d$ is the nominal demand in the corresponding PGLib instance. The proposed convexification procedures are performed in Julia using Gurobi 13.0.1 and executed on a Linux-based high-performance computing cluster with 8 allocated CPU cores and 64 GB of memory. The ACOPF upper bounds and SOCP baseline solutions are computed separately on a MacBook equipped with an Apple M3 chip and 16 GB of memory, using Ipopt 1.14.1 and Mosek 11.0.1, respectively.

For the FBBT procedure, the maximum number of iterations is set to 10, with arctangent envelopes added at the second iteration. The SOC outer-approximation procedure is limited to 20 iterations. The other detailed parameter settings are reported in Table~\ref{tab:offline_parameters}.

\begin{table}[H]
\centering
\caption{Parameter settings for bound tightening and SOC outer approximation.}
\label{tab:offline_parameters}
\small
\setlength{\tabcolsep}{2pt}
\renewcommand{\arraystretch}{0.95}

\begin{tabular}{@{}
    >{\raggedright\arraybackslash}p{0.20\linewidth}
    >{\centering\arraybackslash}p{0.10\linewidth}
    >{\raggedright\arraybackslash}p{0.60\linewidth}
@{}}
\hline
Parameter & Value & Meaning \\
\hline

\multicolumn{3}{c}{\textit{Bound tightening}} \\
\hline

\texttt{max\_iter}
& $10$
& Maximum FBBT iterations. \\

\texttt{bound\_tol}
& $10^{-6}$
& Minimum bound-interval width. \\

\texttt{improve\_tol}
& $10^{-6}$
& Stopping tolerance on $\Delta_k$. \\

\texttt{cut\_tol}
& $10^{-3}$
& Tolerance for degenerate $c$-$s$ boxes. \\

\texttt{cut\_type}
& \texttt{[:AT]}
& Envelope cuts used in FBBT. \\

\texttt{cut\_add\_iter}
& $2$
& Iteration to introduce cut enhancement. \\

\hline
\multicolumn{3}{c}{\textit{SOC outer approximation}} \\
\hline

\texttt{max\_iter}
& $20$
& Maximum SOC-OA iterations. \\

\texttt{cut\_vio}
& $10^{-5}$
& Violation tolerance for cut generation. \\

$p_{\rm thm}$
& $1.00$
& Fraction of violated thermal constraints. \\

$p_{\rm jabr}$
& $0.55$
& Fraction of violated Jabr constraints. \\

\texttt{cut\_age}
& $20$
& Maximum inactive cut age. \\

\texttt{cut\_value}
& $10^{-5}$
& Threshold for dropping inactive cuts. \\

$\varepsilon_{\rm par}$
& $10^{-5}$
& Parallel-cut filtering tolerance. \\

$\varepsilon_{\rm obj}$
& $10^{-5}$
& Objective-improvement tolerance. \\

\hline
\end{tabular}
\end{table}

We use three measures to evaluate the performance of the proposed framework. The \emph{relaxation gap} (Gap), which is the relative difference between an AC-feasible upper bound and the lower bound obtained from the corresponding relaxation, measures relaxation strength after both bound tightening and SOC outer approximation. Computation time (Time) is reported to evaluate their computational cost. For bound tightening, we additionally measure the reduction in the overall $c$-$s$ bounding region. For each line $(i,j)$ at time period $t$, the area of the corresponding bounding box at iteration $k$ is defined as

\begin{equation}
    B^k_{ij,t} = (\bar c^k_{ij,t}-\underline c^k_{ij,t}) (\bar s^k_{ij,t}-\underline s^k_{ij,t}). 
\end{equation}

Let $\mathcal{K}$ denote the set of all $(ij,t)$ pairs. The instance-level area reduction ratio (Area red.) at iteration $k$ is defined as 
\begin{equation}
    \rho^k_{\mathrm{area}} = 1- \frac{ \sum_{(ij,t)\in\mathcal{K}} B^k_{ij,t} }{ \sum_{(i,j,t)\in\mathcal{K}} B^0_{ij,t} }. \label{eq:area_reduction}
\end{equation}

A larger value of $\rho^k_{\mathrm{area}}$ indicates a greater reduction in the overall $c$-$s$ bounding region.

\subsection{Small instance illustration}
\subsubsection{Bounds evolution}
We use the case 3 instance to illustrate the proposed offline convexification and bound tightening framework. The procedure contains two stages. First, the bounds of the $c$-$s$ variables are tightened iteratively. Second, using the tightened bounds with arctangent envelopes enhancement, an offline SOC outer-approximation procedure is applied to generate reusable SOC outer-approximation cuts for constraints~\eqref{eq:thermal_limit} and~\eqref{eq:cs_rotated_cone}.

\begin{figure}[!h]
\centering
\includegraphics[width=0.5\linewidth]{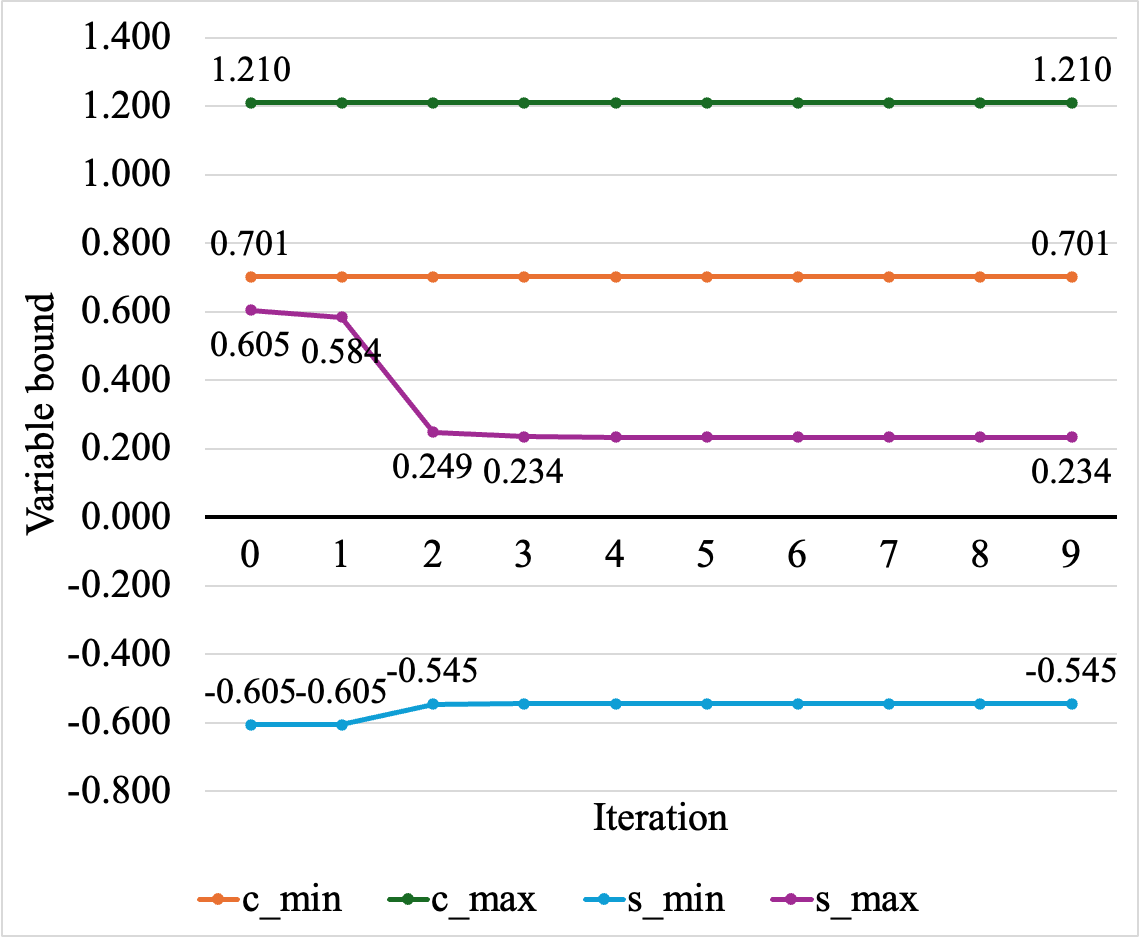}
\caption{Bound evolution for variables $c_{12}$ and $s_{12}$ in the case 3 instance.}
\label{fig:case3}
\end{figure}

Figure~\ref{fig:case3} shows the bound evolution for line $(1,2)$. The initial bounds of $c_{12}$ and $s_{12}$ are $[0.701,1.21]$ and $[-0.605,0.605]$, respectively. After the first bound tightening iteration, the upper bound of $s_{12}$ decreases to $0.584$. Before the second iteration, arctangent cuts are added to strengthen the current relaxation, further tightening the bounds of $s_{12}$ to $[-0.545,0.249]$, after which only marginal improvements are observed. Similar behavior occurs for lines $(1,3)$ and $(2,3)$, and the procedure terminates after 9 iterations because none of the $c$-$s$ variables exhibit a significant improvement beyond the prescribed threshold.

Using the area reduction ratio defined by~\eqref{eq:area_reduction}, the bound region is reduced by $55.26\%$ after the first iteration and $67.59\%$ after the second. This rapid improvement is mainly attributed to the initial bound-tightening step and the addition of arctangent envelopes before the second iteration. After the third iteration, the area reduction ratio stabilizes at around $68.14\%$, indicating that most of the attainable tightening has already been achieved. Meanwhile, the relaxation gap decreases from $8.36\%$ in the first iteration to $7.18\%$ in the final round. 

After obtaining the tightened bounds, we construct the strengthened McCormick-based demand slack model \eqref{model:agg} as the base model. Arctangent envelopes are also included in the relaxation before cut generation. Starting from this enhanced LP relaxation, we perform an offline SOC outer-approximation procedure using Algorithm \ref{alg:offline_soc_cuts}. The procedure generates 11 cuts in total, which are then stored for reuse in the online stage for new instances.

\subsection{Computational results}

\subsubsection{Offline bounds and cuts generation comparison}

We next compare the offline bounds and cuts generated from the standard nominal model and from the proposed demand slack model. The purpose of this comparison is to evaluate whether introducing slack variables preserves the quality of the resulting bounds and cuts while making them valid for a range of demand realizations. We test all nominal PGLib instances with fewer than 1000 buses and the comparison results are reported in Table~\ref{tab:offline_compare}. 

\begin{table}[!h]
\centering
\caption{Comparison of FBBT performance for the nominal and demand slack models.}
\label{tab:offline_compare}
\small
\setlength{\tabcolsep}{5.0pt}
\renewcommand{\arraystretch}{1.05}
\begin{tabular}{l|rrr|rrr}
\hline
& \multicolumn{3}{c|}{BT on nominal model}
& \multicolumn{3}{c}{BT on demand slack model} \\
Instance
& Iter.
& Area red.
& Time (s)
& Iter.
& Area red.
& Time (s) \\
\hline
3-lmbd       & 10 & 83.60\% & 7.28     & 9  & 68.10\% & 12.23 \\
5-pjm        & 10 & 96.30\% & 0.99     & 10 & 94.50\% & 0.97 \\
14-ieee      & 10 & 88.10\% & 2.62     & 10 & 85.80\% & 2.79 \\
24-rts       & 10 & 89.70\% & 7.85     & 10 & 84.40\% & 8.07 \\
30-as        & 10 & 94.50\% & 11.48    & 10 & 95.40\% & 12.18 \\
30-ieee      & 10 & 90.30\% & 10.73    & 10 & 90.80\% & 11.81 \\
39-epri      & 10 & 94.70\% & 12.25    & 10 & 90.70\% & 13.45 \\
57-ieee      & 10 & 81.20\% & 74.84    & 10 & 81.10\% & 82.82 \\
60-c         & 10 & 81.80\% & 39.56    & 10 & 77.70\% & 38.35 \\
73-rts       & 10 & 85.20\% & 134.17   & 10 & 80.40\% & 143.26 \\
89-pegase    & 10 & 69.80\% & 1785.50  & 10 & 68.20\% & 1894.71 \\
118-ieee     & 10 & 85.70\% & 688.19   & 10 & 84.30\% & 725.20 \\
162-dtc      & 10 & 85.80\% & 5800.92  & 10 & 83.70\% & 6275.34 \\
179-goc      & 10 & 79.40\% & 837.12   & 10 & 79.20\% & 870.18 \\
197-snem     & 10 & 79.30\% & 448.10   & 10 & 92.40\% & 574.44 \\
200-activ    & 10 & 77.70\% & 1606.87  & 10 & 84.80\% & 1793.38 \\
240-pserc    & 10 & 71.70\% & 4495.09  & 10 & 70.80\% & 4523.06 \\
300-ieee     & 10 & 79.30\% & 6634.60  & 10 & 79.30\% & 6967.09 \\
500-goc      & 10 & 85.90\% & 29420.74 & 10 & 86.50\% & 31587.30 \\
588-sdet*     & 1  & 0.00\%  & 67.94    & 10 & 90.30\% & 11032.04 \\
793-goc*      & 2  & 81.90\% & 198.12   & 10 & 91.10\% & 29900.13 \\
\hline
Average      & -- & 80.09\% & 2489.76
             & -- & 83.78\% & 4593.75 \\
Average (excl. *)      & -- & 84.21\% & 2489.76 & -- & 83.06\% & 2922.98 \\
\hline
\end{tabular}

\vspace{1mm}
\begin{flushleft}
\footnotesize
\emph{Note:} ``BT'' denotes FBBT. Area reductions are reported as
percentages, and time denotes the total bound-tightening runtime. "*" denotes instances for which  numerical issues were encountered during the FBBT process. The ``Average (excl. *)'' row excludes \texttt{588-sdet} and
\texttt{793-goc}.
\end{flushleft}
\end{table}

As shown in Table~\ref{tab:offline_compare}, we first apply the convexification framework to the standard nominal model without slack variables. On average, the procedure reduces the bounding-box area of the $c$--$s$ variables by $80.09\%$ and requires $2489.76$ seconds. For most test cases, FBBT reaches the prescribed maximum of ten iterations. However, numerical difficulties cause the procedure to terminate after only one iteration for \texttt{case588-sdet} and two iterations for \texttt{case793-goc}.

We next apply the same FBBT procedure to the demand slack model. The average area reduction increases from $80.09\%$ for the nominal model to $83.78\%$ for the demand slack model, while the average computation time increases from $2489.76$ seconds to $4593.75$ seconds. For most individual instances, the demand slack model produces a slightly smaller area reduction than the nominal model, as expected, because allowing the demand to
vary enlarges the feasible region. Nevertheless, the resulting reductions remain close for many instances, indicating that the feasible ranges of the $c$--$s$ variables do not change substantially over the considered demand range.

The larger average area reduction obtained from the demand slack model is mainly driven by several numerically challenging instances. In particular, for \texttt{case197-snem}, \texttt{case200-activ}, \texttt{case588-sdet}, and \texttt{case793-goc}, the demand slack model achieves substantially greater area reduction than the nominal model. This result is caused by the
numerical instability of FBBT on the nominal model. For these instances, some variables obtain nearly identical lower and upper bounds, resulting in near-degenerate intervals. To preserve numerical stability, a bound is not updated when the resulting interval width falls below the tolerance $10^{-6}$. In the most severe cases, \texttt{case588-sdet} and \texttt{case793-goc}, these numerical issues cause the nominal-model FBBT procedure to terminate prematurely, yielding area reductions of only $0.00\%$ and $81.90\%$, respectively. In contrast, the slightly larger feasible region of the demand slack model avoids many of these near-degenerate bounds, allowing FBBT to complete all ten iterations and achieve area reductions of $90.30\%$ and $91.10\%$, respectively. This also shows the numerical robustness of the demand slack model.

\subsubsection{Comparison between offline and online convexification}

\begin{table*}[!h]
\centering
\caption{Comparison of online and offline bound-tightening and SOC-OA strategies}
\label{tab:online_offline_compare}
\small
\setlength{\tabcolsep}{2.5pt}
\renewcommand{\arraystretch}{1.08}

\resizebox{\textwidth}{!}{
\begin{tabular}{l|rrr|rrr|rrr}
\hline
& \multicolumn{3}{c|}{Online FBBT + online SOC-OA}
& \multicolumn{3}{c|}{Offline bounds + offline cuts}
& \multicolumn{3}{c}{Offline bounds + online SOC-OA} \\
Instance
& Gap (\%)
& Mean $\pm$ Std. gap (\%)
& Time (s)
& Gap (\%)
& Mean $\pm$ Std. gap (\%)
& Time (s)
& Gap (\%)
& Mean $\pm$ Std. gap (\%)
& Time (s) \\
\hline

case3\_lmbd
& 1.47 & $1.54 \pm 0.61$ & 1.72
& 3.05 & $3.13 \pm 0.82$ & 0.26
& 2.90 & $2.97 \pm 0.70$ & 0.34 \\

case5\_pjm
& 0.25 & $0.34 \pm 0.59$ & 1.27
& 0.38 & $0.72 \pm 1.72$ & 0.04
& 0.36 & $0.67 \pm 1.58$ & 0.06 \\

case14\_ieee
& 3.33 & $3.47 \pm 1.12$ & 2.81
& 3.42 & $3.57 \pm 1.15$ & 0.11
& 3.35 & $3.49 \pm 1.13$ & 0.12 \\

case24\_ieee\_rts
& 1.69 & $1.70 \pm 0.14$ & 9.22
& 3.16 & $3.18 \pm 0.41$ & 0.37
& 2.44 & $2.45 \pm 0.25$ & 0.79 \\

case30\_as
& 16.28 & $18.21 \pm 9.87$ & 13.05
& 19.70 & $21.60 \pm 10.51$ & 0.36
& 17.72 & $19.59 \pm 10.05$ & 0.81 \\

case30\_ieee
& 1.61 & $1.61 \pm 0.10$ & 11.26
& 1.93 & $1.93 \pm 0.10$ & 0.25
& 1.81 & $1.81 \pm 0.10$ & 0.29 \\

case39\_epri
& 0.44 & $0.46 \pm 0.27$ & 13.41
& 0.64 & $0.66 \pm 0.29$ & 0.28
& 0.54 & $0.57 \pm 0.28$ & 0.33 \\

case57\_ieee
& 3.64 & $3.64 \pm 0.09$ & 76.93
& 3.86 & $3.86 \pm 0.12$ & 0.78
& 3.82 & $3.82 \pm 0.12$ & 0.88 \\

case60\_c
& 1.18 & $1.21 \pm 0.36$ & 40.08
& 1.70 & $1.77 \pm 0.62$ & 0.76
& 1.37 & $1.44 \pm 0.57$ & 0.85 \\

case73\_ieee\_rts
& 0.96 & $0.99 \pm 0.41$ & 145.77
& 1.81 & $1.86 \pm 0.57$ & 2.45
& 1.33 & $1.36 \pm 0.43$ & 4.96 \\

case89\_pegase
& 3.86 & $6.80 \pm 7.19$ & 1835.00
& 4.12 & $6.98 \pm 7.18$ & 9.75
& 3.96 & $6.88 \pm 7.19$ & 11.01 \\

case118\_ieee
& 15.75 & $15.77 \pm 1.00$ & 710.48
& 17.68 & $17.70 \pm 0.89$ & 5.07
& 17.16 & $17.19 \pm 0.92$ & 5.14 \\

case162\_ieee\_dtc
& 4.08 & $4.69 \pm 2.18$ & 5810.24
& 5.21 & $5.63 \pm 2.09$ & 15.95
& 4.14 & $4.76 \pm 2.22$ & 17.31 \\

case179\_goc
& 2.97 & $3.01 \pm 0.59$ & 831.64
& 3.92 & $3.97 \pm 0.76$ & 4.98
& 2.66 & $2.79 \pm 0.97$ & 5.52 \\

case197\_snem
& 0.15 & $0.16 \pm 0.21$ & 454.71
& 0.42 & $0.43 \pm 0.05$ & 5.11
& 0.37 & $0.37 \pm 0.05$ & 5.55 \\

case200\_activ
& 0.06 & $0.06 \pm 0.00$ & 1614.84
& 0.07 & $0.07 \pm 0.04$ & 9.44
& 0.06 & $0.06 \pm 0.00$ & 54.77 \\

case240\_pserc
& 0.80 & $0.82 \pm 0.31$ & 4495.49
& 0.92 & $0.94 \pm 0.33$ & 13.79
& 0.82 & $0.84 \pm 0.31$ & 15.14 \\

case300\_ieee
& 0.72 & $0.73 \pm 0.22$ & 6477.23
& 1.12 & $1.15 \pm 0.35$ & 12.88
& 0.81 & $0.82 \pm 0.26$ & 13.10 \\

case500\_goc
& 2.74 & $2.89 \pm 1.28$ & 30648.99
& 2.93 & $3.08 \pm 1.29$ & 23.53
& 2.77 & $2.92 \pm 1.28$ & 35.63 \\

case588\_sdet
& 2.80 & $2.81 \pm 0.25$ & 115.01
& 2.12 & $2.12 \pm 0.08$ & 38.54
& 2.05 & $2.05 \pm 0.08$ & 39.80 \\

case793\_goc
& 8.50 & $8.94 \pm 3.50$ & 584.74
& 6.39 & $6.42 \pm 0.71$ & 53.05
& 6.28 & $6.30 \pm 0.67$ & 283.53 \\

\hline
\textbf{Average}
& \textbf{3.49}
& $\mathbf{3.80 \pm 1.44}$
& \textbf{2566.38}
& \textbf{4.03}
& $\mathbf{4.32 \pm 1.43}$
& \textbf{9.42}
& \textbf{3.65}
& $\mathbf{3.96 \pm 1.39}$
& \textbf{23.62} \\
\hline

\end{tabular}
}
\begin{flushleft}
\footnotesize
\emph{Note:} Gap and runtime are reported using shifted geometric means (SGMs), with shifts of one percentage point and one second, respectively. Mean $\pm$ Std. reports the arithmetic mean and standard deviation of the gap over the sampled instances.
\end{flushleft}
\end{table*}

\begin{figure}[!h]
    \centering
    \includegraphics[width=0.5\linewidth]{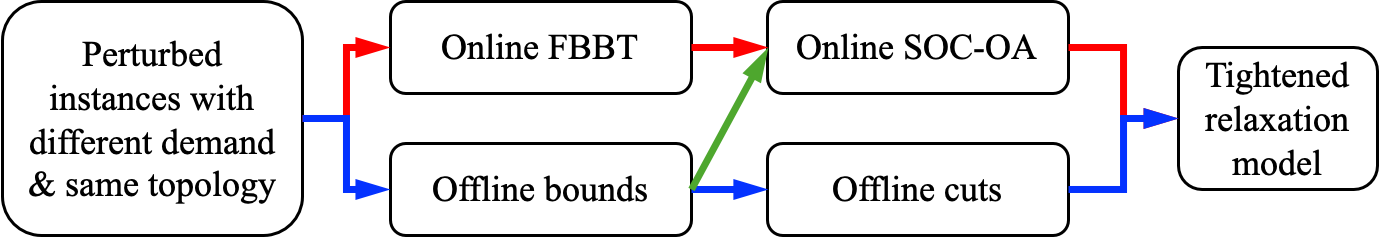}
    \caption{Different experiment strategies for offline and online comparison}
    \label{fig:exp}
\end{figure}

We further evaluate the effectiveness of offline bounds and cuts on randomly perturbed demand instances. For each PGLib case, 10 demand profiles are randomly selected from the 100 generated load samples. Three strategies are compared, as illustrated in Figure~\ref{fig:exp}. 
The first strategy, \emph{Online BT + Online SOC-OA}, applies the complete strengthen framework independently to each perturbed instance, including both bound tightening and SOC outer approximation, without introducing slack variables. The second strategy, \emph{Offline Bounds + Offline Cuts}, directly applies the bounds and cuts precomputed from the demand slack model to the McCormick relaxation with arctangent envelopes. The third strategy, \emph{Offline Bounds + Online SOC-OA}, reuses the offline precomputed bounds but regenerates the SOC outer-approximation cuts online for each perturbed instance. 

The relaxation gap and computation time for the three strategies are reported in Table~\ref{tab:online_offline_compare}. To summarize relaxation gaps and computation times across demand profiles, we use the shifted geometric mean (SGM) defined for nonnegative observations $z_1,\ldots,z_N$ as

\[\operatorname{SGM}_s(z_1,\ldots,z_N) = \left(\prod_{r=1}^{N}(z_r+s)\right)^{1/N}-s,\] 
where $s>0$ is a prescribed shift. We use a shift of $1\%$ for the relaxation gap and 1 second for computation time. For each case, the SGM is computed over the 10 sampled demand profiles. To additionally characterize the variation in relaxation quality, Table~\ref{tab:online_offline_compare} also reports the arithmetic mean and standard deviation of the gap. 

The results show that the offline bounds recover most of the relaxation strength obtained by online FBBT. The fully online strategy provides the strongest overall relaxation, with an average case-level SGM gap as $3.49\%$ and an arithmetic mean gap of $3.80\%$, with an average standard deviation of $1.44\%$. However, because FBBT and SOC-OA are performed independently for each perturbed demand profile, its case-average $\mathrm{SGM}$ computational time is $2566.38$ seconds.

In comparison, directly applying the offline bounds and offline cuts yields a $\mathrm{SGM}$ gap of $4.03\%$ and an arithmetic mean gap of $4.32\%$ with a standard deviation of $1.43\%$. Although the relaxation is slightly weaker than that obtained by the fulling online strategy, the case-shifted $\mathrm{SGM}$ time is only $9.42$ seconds. Thus, the online computational time is reduced by approximately $272\times$, while the arithmetic mean gap increases by only $0.52$ percentage points. This demonstrates that the precomputed bounds and cuts retain most of the relaxation strength of the fully online procedure while substantially reducing the cost of repeated scenario evaluations.

The hybrid strategy, which uses offline bounds but regenerates SOC-OA cuts online, provides a useful compromise between relaxation strength and computation speed. It obtains a $\mathrm{SGM}$ gap of $3.65\%$ and an arithmetic mean gap of $3.96\%$ with a standard deviation of $1.39\%$. These values are close to the corresponding fully online gaps of $3.49\%$ and $3.80\%$, respectively. At the same time, its case-average $\mathrm{SGM}$ computational time is only $23.62$ seconds, which is still more than $108\times$ faster than the fully online strategy. Therefore, most of the relaxation strength can be recovered by reusing the offline bounds and performing only the relatively inexpensive SOC-OA procedure online.

\begin{figure}[!h]
    \centering
    \includegraphics[width=0.5\linewidth]{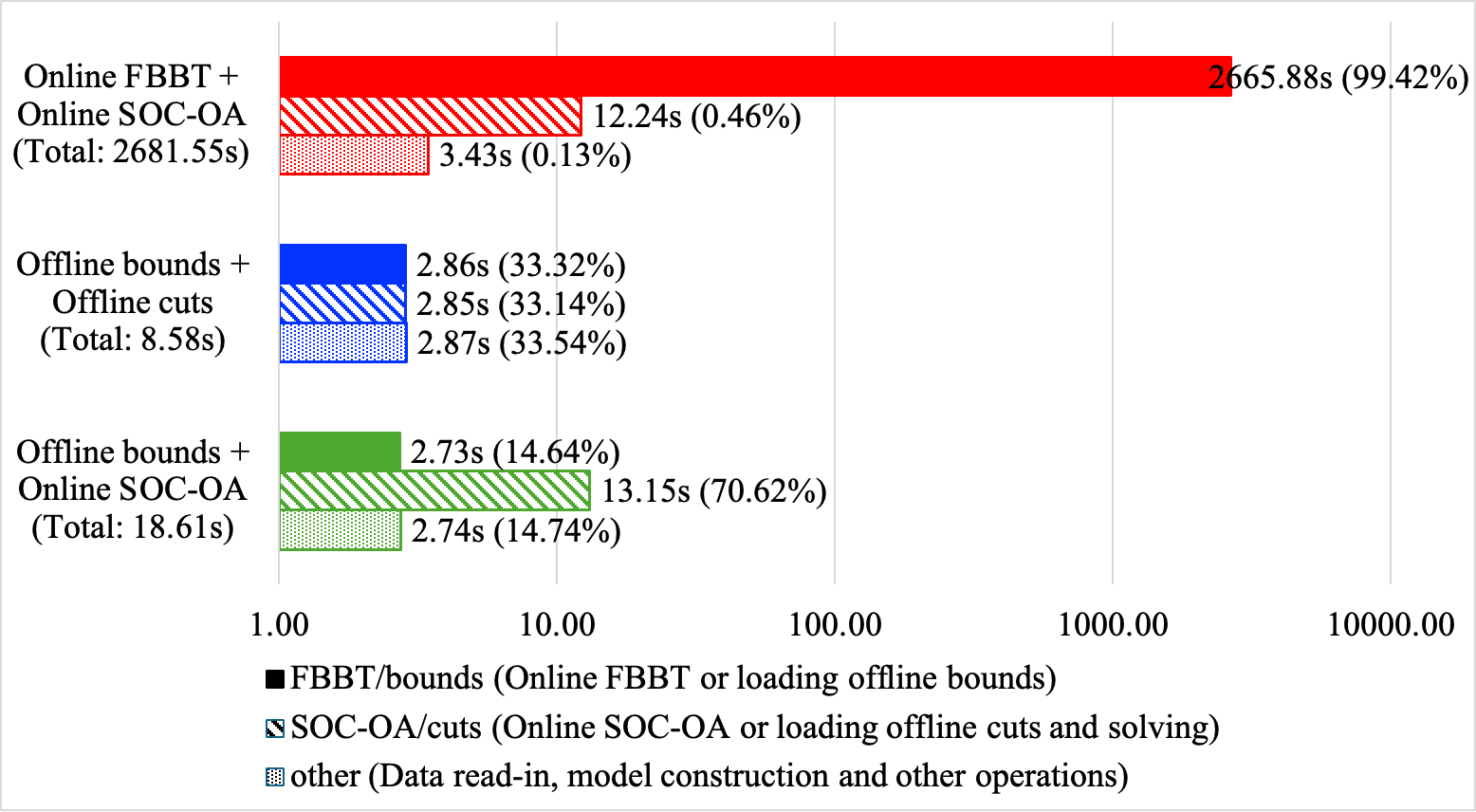}
    \caption{Percentage breakdown of computational time across different strategies}
    \label{fig:time_portion}
\end{figure}

To further examine the source of computational savings, Figure~\ref{fig:time_portion} reports the component-wise average runtime breakdown of the three strategies on a logarithmic scale. Unlike the SGM times reported in Table~\ref{tab:online_offline_compare}, arithmetic mean times are used here because the individual runtime components should be additive and can therefore be consistently decomposed. For each strategy, the total runtime is divided into 3 parts: \emph{FBBT/bounds} time includes either the complete online FBBT procedure or, for offline strategies, the time required to build the relaxation, load and apply those bounds for the tightened relaxation. The \emph{SOC-OA/cuts} component includes either the complete online SOC-OA time required to load and apply the precomputed cuts and solve the resulting relaxation. The remaining component includes data read in, and other essential operations that are not separately timed. The total average runtime is reported alongside each strategy, while each component is labeled by both its average time and its percentage of the corresponding total.

The fully online strategy requires an average of $2681.55$ seconds, of which FBBT accounts for $2665.88$ seconds, or $99.42\%$. In comparison, online SOC-OA requires only $12.24$ seconds ($0.46\%$), while model construction, data handling, and the remaining operations require $3.43$ seconds($0.13\%$). This large difference confirms that repeatedly performing instance-specific FBBT is the primary computational bottleneck of the fully online procedure.

Reusing the offline bounds removes most of this cost. For the \emph{Offline Bounds + Offline Cuts} strategy, the total average runtime decreases to $8.58$ seconds. Loading and applying the offline bounds, processing the offline cuts and solving the resulting relaxation, and the remaining operations require $2.86$, $2.85$, and $2.87$ seconds, respectively, corresponding to $33.32\%$, $33.14\%$, and $33.54\%$ of the total runtime. 
For the hybrid \emph{Offline Bounds + Online SOC-OA} strategy, the total average runtime is $18.61$ seconds. Applying the offline bounds requires only $2.73$ seconds ($14.64\%$), while regenerating the SOC-OA cuts online becomes the dominant component at $13.15$ seconds ($70.62\%$); the remaining operations require $2.74$ seconds ($14.74\%$). Despite performing SOC-OA online, the hybrid strategy remains substantially less expensive than the fully online procedure because the costly instance-specific bound-tightening stage is avoided.

Overall, Table~\ref{tab:online_offline_compare} shows that the offline framework substantially reduces the online computational burden by eliminating repeated FBBT while preserving most of the relaxation quality. The fully offline setting is the fastest option, while the hybrid setting provides a stronger relaxation with only a modest increase in online time. Therefore, the offline bounds and cuts are particularly useful when many demand scenarios must be evaluated repeatedly.

\section{Conclusion}

This paper proposes an offline convexification and bound tightening framework for strengthening linear relaxations of the ACOPF problem. The framework combines bound tightening, arctangent-envelope enhancement, and SOC outer-approximation cut generation. The main motivation is that, in many operational settings, a family of ACOPF instances shares the same network topology and differs primarily in load profiles. Instead of performing expensive bound tightening independently for each new instance, the proposed approach constructs a demand slack model over prescribed load ranges and computes tightened bounds and valid cuts offline. These bounds and cuts can then be reused for future scenarios whose demands fall within the specified ranges.

Computational results on PGLib instances show that 
% the proposed framework significantly improves relaxation quality. 
by using the demand slack model, the tightened bounds are representative of the feasible $c$-$s$ region across perturbed demand instances. Reusing the offline bounds and cuts produces an average gap of $4.03\%$, which is close to the $3.49\%$ gap obtained by the fully online convexification procedure, while achieving a $272\times$ speedup. When the offline bounds are combined with online SOC outer approximation, the final average gap is further reduced to $3.65\%$.

Overall, the results indicate that the proposed offline convexification and bound tightening framework provides a practical way to balance the cost of bound tightening and cut generation across related ACOPF instances. The offline bounds capture most of the gap-closing effect of online FBBT, while the offline cuts further strengthen the resulting LP relaxation. Since the generated bounds and cuts are valid over prescribed demand ranges, they can be directly reused for repeated scenario evaluations. This makes the framework well suited for multi-period planning, contingency analysis, and other applications where many ACOPF instances with the same network structure must be solved repeatedly.

\bibliographystyle{IEEEtran}
\bibliography{FBBT_refs}
% Insert the remaining paper here, starting with the Introduction.
\end{document}